\documentclass{gtart}     
\usepackage{amssymb,amsmath}
\newtheorem{thm}{Theorem}[section]
\newtheorem{conj}[thm]{Conjecture}

\newtheorem{prop}[thm]{Proposition}
\newtheorem{cor}[thm]{Corollary}
\newtheorem{rem}[thm]{Remark}

\newcommand{\bC}{\mathbb C}

\newcommand{\blowup}{\overline{\bC P^2}}

\DeclareMathOperator{\CPB}{\overline{\C P}}

\newcommand{\dvolg}{{d\text{vol}_g}}

\newcommand{\C}{\mathbb{ C}}

\newcommand{\R}{\mathbb{ R}}

\newcommand{\Z}{\mathbb{ Z}}

\newcommand{\Spc}{\mathrm{Spin}^c}

\newcommand{\cH}{{\mathcal H}}

\newcommand{\FAP}{{F_{\hat{A}}^+}}
\newcommand{\FAM}{{F_{\hat{A}}^-}}

\newcommand{\PT}{{\widetilde{P}}}

\begin{document}
%
\input gtoutput
\volumenumber{2}\papernumber{1}\volumeyear{1998}
\pagenumbers{1}{10}\published{16 January 1998}
\proposed{Peter Kronheimer}\seconded{Ronald Stern, Gang Tian}
\received{8 September 1997}\revised{14 January 1998}
\accepted{15 January 1998}
\title{Einstein metrics and smooth structures}         

\author{D Kotschick}

\address{Mathematisches Institut\\Universit\"at Basel\\Rheinsprung 21\\%
4051 Basel, Switzerland}

\asciiaddress{Mathematisches Institut, Universitat Basel, Rheinsprung 21,
4051 Basel, Switzerland}

\email{dieter@math-lab.unibas.ch}

\begin{abstract} 

We prove that there are infinitely many pairs of 
homeomorphic non-diffeo\-morphic smooth $4$--manifolds,
such that in each pair one manifold admits an Einstein
metric and the other does not. We also show that there
are closed $4$--manifolds with two smooth structures
which admit Einstein metrics with opposite signs of
the scalar curvature.

\end{abstract}

\asciiabstract{We prove that there are infinitely many pairs of 
homeomorphic non-diffeomorphic smooth 4-manifolds,
such that in each pair one manifold admits an Einstein
metric and the other does not. We also show that there
are closed 4-manifolds with two smooth structures
which admit Einstein metrics with opposite signs of
the scalar curvature.}



\maketitlepage

%

In dimensions strictly smaller than four Einstein metrics
have constant curvature and are therefore rare.
In dimension four Einstein metrics of non-constant
curvature exist, but it is still the case that existence of 
such a metric imposes non-trivial restrictions on the 
underlying manifold\footnote{No such restrictions are known
in higher dimensions.}.
For closed orientable Einstein $4$--manifolds $X$
the Euler characteristic has to be non-negative,
and, furthermore, the Hitchin--Thorpe inequality
\begin{equation}\label{eq:HT}
e (X) \geq \frac{3}{2}\vert\sigma (X)\vert 
\end{equation}
must hold, where $e$ denotes the Euler characteristic
and $\sigma$ the signature. This condition is very
crude, and is certainly homotopy invariant, as
are the restrictions coming from Gromov's
notion of simplicial volume~\cite{gromov}, and from 
the existence of maps of non-zero degree to hyperbolic 
manifolds~\cite{BCG}.

Our aim in this note is to discuss existence and 
non-existence of Einstein metrics as a property of 
the smooth structure. We shall exhibit infinitely many pairs of 
homeomorphic non-diffeomorphic smooth $4$--manifolds,
such that in each pair one manifold admits an Einstein
metric and the other does not. This shows for the 
first time that the smooth structures of $4$--manifolds
form definite obstructions to the existence of an
Einstein metric.

An isolated example of such a pair can be obtained as follows. 
Hitchin~\cite{HT} showed that Einstein manifolds for 
which~\eqref{eq:HT} is an equality are either flat or quotients 
of a $K3$ surface with a Calabi--Yau metric. Thus, the 
existence of smooth manifolds homeomorphic but not
diffeomorphic to the $K3$ surface, which is known from Donaldson 
theory~\cite{FM} and also follows easily from Seiberg--Witten 
theory, see for example~\cite{bams,KKM}, shows that by changing only
the differentiable structure one can pass from a manifold with 
an Einstein metric to one without. The point of our examples
is that there are lots of them, and they do not arise
from the borderline case of a non-existence result.
They are in some sense generic.

We shall also discuss a conjecture concerning uniqueness
of Einstein metrics on $4$--manifolds which complements
the discussion of existence. This too depends on a 
consideration of different smooth structures on a fixed
topological manifold.

\section{Smooth structures as obstructions}

We shall use Seiberg--Witten invariants to show that certain
smooth structures obstruct the existence of Einstein metrics,
and refer the reader to~\cite{witten,bams,KKM} for the 
definitions and basic properties of the invariants.
All manifolds in this section are closed, smooth, oriented 
$4$--manifolds. For the sake of simplicity, we assume $b_2^+>1$ 
throughout, though this is not essential.
 
We shall need the following result concerning the behaviour
of the invariants under connected summing with $\blowup$.
This is usually referred to as a blowup formula.

\begin{prop}[\cite{FS,KMT}]\label{p:blowup}
Let $\PT (Y)$ be a $\Spc$--structure on $Y$, and  $X = Y \# \CPB^2$,
with $E$ a generator of $H^2(\CPB^2,\Z )$. Then $X$ has a 
$\Spc$--structure $\PT (X)$ with $c_1(\PT (X)) = c_1(\PT (Y)) + E$, 
such that the Seiberg--Witten invariants of $\PT (Y)$ and of $\PT (X)$
are equal (up to sign).
\end{prop}
As the reflection in $E^{\perp}$ in the cohomology of $X$ is 
realised by a self-diffeo\-morphism, the naturality of the 
invariants shows that there is another $\Spc$--structure with the
same Seiberg--Witten  invariant, up to sign, and with\break
$c_1(\PT (X)) = c_1(\PT (Y)) - E$.

Using this, we can prove the following version of a 
theorem of LeBrun~\cite{lebrun}:

\begin{thm}\label{t:blowup}
Let $Y$ be a manifold with a non-zero Seiberg--Witten invariant
(of any degree), and 
$X = Y \# k \CPB^2$. If $k > \frac{2}{3}(2e(Y) + 3\sigma(Y))$, then
$X$ does not admit an Einstein metric.
\end{thm}
\begin{proof}
If $\PT (X)$ has a non-zero Seiberg--Witten invariant, 
then for every Riemannian metric $g$ there must be a solution
$(A,\phi )$ of the monopole equations. Denoting by $\hat A$ the
connection induced by the $\Spc$--connection $A$ on the 
determinant bundle of the spinor bundle, we  have
\begin{align*} 
c_1^2(\PT (X)) 
&= \frac{1}{4\pi^2} \int_X (|\FAP|^2 - |\FAM|^2) \dvolg
\le \frac{1}{4\pi^2} \int_X |\FAP|^2  \dvolg  \\
&=\frac{1}{32\pi^2}\int_X\vert\phi\vert^4\dvolg
\le \frac{1}{32\pi^2} \int_X s_g^2 \dvolg \ , 
\end{align*}
where $s_g$ denotes the scalar curvature of $g$. 

Given any class $c \in H^2(X,\R )$, denote by $c^+$ the projection 
of $c$ into the subspace $\cH^2_+ \subset \cH^2$ of $g$--self-dual 
harmonic forms along the subspace $\cH^2_-$ of $g$--anti-self-dual 
harmonic forms. The argument above really proves 
\[ \left( c_1(\PT (X))^+ \right)^2 \le  \frac{1}{32\pi^2} \int_X s_g^2 
\dvolg \ . \]

If $\PT (Y)$ is a $\Spc$--structure on $Y$ with non-zero
Seiberg--Witten invariant, then, by Proposition~\ref{p:blowup}
and the subsequent remark, there are $\Spc$--structures $\PT (X)$ on
$X = Y \# k \CPB^2$ with non-zero Seiberg--Witten invariants
and with 
\[ c_1(\PT (X)) = c_1(\PT (Y)) + \sum_{i=1}^k (-1)^{\epsilon_i}E_i \]
for any choice of the signs $(-1)^{\epsilon_i}$.
Choose the signs so that 
$$
(-1)^{\epsilon_i}E_i^+\cdot c_1(\PT (Y))^+\ge 0 \ .
$$
Then
\begin{multline*}
\frac{1}{32\pi^2} \int_X s_g^2 \dvolg  
\ge  \left( c_1(\PT (X))^+ \right)^2 \\
= \left( c_1(\PT (Y))^+ \right)^2 +     
2\sum_{i=1}^k (-1)^{\epsilon_i}E_i^+ \cdot c_1(\PT (Y))^+
+\left( \sum_{i=1}^k (-1)^{\epsilon_i}E_i^+ \right)^2 \\
\ge  \left( c_1(\PT (Y))^+ \right)^2 \ge  c_1(\PT (Y))^2                                         \ge 2e(Y) + 3\sigma(Y)  \\                                
= 2(e(X) - k) + 3(\sigma(X) + k)                             
= 2e(X) + 3\sigma(X) + k \ ,
\end{multline*}
where we have used the inequality $c_1(\PT (Y))^2\ge 2e(Y) +3\sigma(Y)$
which is equivalent to the assertion that the moduli space associated
with $\PT (Y)$ has non-negative dimension.

Thus, we have proved $\frac{1}{32\pi^2}\int_X s_g^2 \dvolg\ge 
2e(X)+3\sigma(X)+k$ for every metric $g$ on $X$.

Suppose now that $g$ is Einstein. Then the Chern--Weil integrals
for the Euler characteristic and the signature of $X$ give
\begin{align*}
2e(X) + 3\sigma(X)  
&= \frac{1}{4\pi^2} \int_X (\frac{1}{24}s_g^2 + 2|W_+|^2) \dvolg \\
&\ge \frac{1}{96\pi^2} \int_X s_g^2 \dvolg                         \\
&\ge \frac{1}{3} ( 2e(X) + 3\sigma(X) + k ) \ ,
\end{align*}
where $W_+$ denotes the self-dual part of the Weyl tensor
of $g$. Therefore $k\le 2( 2e(X)+3\sigma(X))=2(2e(Y)+3\sigma(Y)-k)$,
which implies $k \le \frac{2}{3} ( 2e(Y) + 3\sigma(Y))$.
\end{proof}

Theorem~\ref{t:blowup} was proved by LeBrun~\cite{lebrun}, who
also discussed the borderline case $k=\frac{2}{3} (2e(Y)+3\sigma(Y))$,
in the case where $Y$ is complex or symplectic. In that 
case the blown up manifold $X$ is also complex, respectively
symplectic, so that Proposition~\ref{p:blowup} is not
needed.

The following is the main result of this section, giving the
examples mentioned in the introduction.

\begin{thm}\label{t:main}
There are infinitely many pairs $(X_i,Z_i)$ of simply connected
closed oriented smooth $4$--manifolds such that:
\begin{itemize}
\item[\bf1\rm)] $X_i$ is homeomorphic to $Z_i$,
\item[\bf2\rm)] if $i\neq j$, then $X_i$ and $X_j$ are not homotopy
equivalent,
\item[\bf3\rm)] $Z_i$ admits an Einstein metric but $X_i$ does not,
\item[\bf4\rm)] $e (X_i) > \frac{3}{2}\vert\sigma (X_i)\vert$.
\end{itemize}
\end{thm}
\noindent
Note that 3) implies in particular that $X_i$ and $Z_i$
are not diffeomorphic. 
\begin{proof}
We claim that there are simply connected minimal complex 
surfaces $Y_i$, $Z_i$ of general type such that if we  
take $X_i=Y_i\# k\blowup$, for a suitable $k$ with 
$k>\frac{2}{3}(2e(Y_i)+3\sigma(Y_i))$, then the pairs $(X_i,Z_i)$ 
have all the desired properties. The last property, the strict
Hitchin--Thorpe inequality, follows from the Noether
and Miyaoka--Yau inequalities for $Z_i$, which, by the 
first property, has the same Euler characteristic 
and signature as $X_i$.

If we take $Z_i$ to have ample canonical bundle, then
the results of Aubin and Yau on the Calabi conjecture 
show that $Z_i$ admits a K\"ahler--Einstein metric,
compare~\cite{einstein}.
On the other hand, $X_i$ does not admit any Einstein
metric by Theorem~\ref{t:blowup}.

The crucial issue then is to arrange that $Z_i$, with 
ample canonical bundle, is  homeomorphic to the $k$--fold
blowup of $Y_i$, with  $k>\frac{2}{3}(2e(Y_i)
+3\sigma(Y_i))$. As $X_i$ will be automatically non-spin,
$X_i$ and $Z_i$ will be homeomorphic by Freedman's 
classification~\cite{freed} as soon as $Z_i$ is non-spin 
and has the same Euler characteristic and the same signature
as $X_i$. One can find suitable surfaces 
using the known results on the geography of surfaces
of general type, see~\cite{orient} for a summary of
the results.

To exhibit concrete examples, instead of working with the 
topological Euler characteristic and the signature, we shall 
use the first Chern number $c_1^2=2e+3\sigma$
and the Euler characteristic of the structure sheaf
$\chi =\frac{1}{4}(e+\sigma )$.

Under blowing up, $c_1^2$ drops by one and $\chi$ is 
constant. Thus, the Miyaoka--Yau inequality for $Y_i$ 
implies $c_1^2(X_i)<3\chi(X_i)$. In fact, $c_1^2(X_i)$
will be smaller still, because simply connected surfaces
$Y_i$ are not known to exist if we get too close to
the Miyaoka--Yau line $c_1^2=9\chi$. 

The minimal surface $Z_i$ satisfies the same inequalities 
on its characteristic numbers as $X_i$. If the canonical
bundle of $Z_i$ is very ample, then Castelnuovo's theorem,
see~\cite{BPV} page 228, gives $c_1^2(Z_i)\ge 3\chi(Z_i)-10$,
which will contradict the above upper bound for $c_1^2(X_i)$.
Thus, $Z_i$ must be chosen to have ample but not very
ample canonical bundle, and will be in the sector
where 
$$
2\chi(Z_i)-6\leq c_1^2(Z_i) < 3\chi(Z_i) \ ,
$$
the first being the Noether inequality.
The results of Xiao Gang and Z~Chen, cf~\cite{orient},
show that all non-spin simply connected surfaces $Z_i$ with 
ample canonical bundle which are in this sector, and not too
close to the line $c_1^2=3\chi$, will have companions
$X_i$ as required, obtained by blowing up minimal surfaces
$Y_i$. Note that by Beauville's theorem on the canonical
map, cf~\cite{BPV} page 228, all the $Z_i$ will be 
double covers of ruled surfaces.

We can avoid using the results of Xiao and Chen by 
taking for $Z_i$ the following family of Horikawa surfaces,
cf~\cite{BPV}. Let $\Sigma_i$ be the Hirzebruch surface
whose section at infinity $S$ has self-intersection $-i$, 
and let $Z_i$ be a double cover of $\Sigma_i$ branched
in a smooth curve homologous to $B=6S+2(2i+3)F$, where $F$
is the class of the fiber. The double cover is simply 
connected as $B$ is ample, and $K_{Z_i}=\pi^*(K_{\Sigma_i}+
\frac{1}{2}B)=\pi^*(S+(i+1)F)$ is not $2$--divisible
and so $Z_i$ is not spin. Moreover, $K_{Z_i}$ is the 
pullback of an ample line bundle and therefore ample,
so that $Z_i$ admits a K\"ahler--Einstein metric.
The characteristic numbers of $Z_i$ are $c_1^2(Z_i)=
2i+4$ and $\chi(Z_i)=i+5$.

Now, by the classical geography results of Persson~\cite{per},
for all $i$ large enough there are simply connected surfaces 
$Y_i$ of general type with $c_1^2(Y_i)=6i+13$ and $\chi (Y_i)
=i+5$, so that the $(4i+9)$--fold blowup $X_i$ of $Y_i$
is homeomorphic to $Z_i$. 

The pairs $(X_i,Z_i)$ have all the desired properties.
\end{proof}

\begin{rem}\rm
The examples of manifolds without Einstein metrics given by
LeBrun~\cite{lebrun}, namely blowups of hypersurfaces in 
$\C P^3$, cannot be used to prove Theorem~\ref{t:main} because 
they violate the Noether inequality. They are therefore not 
homeomorphic to minimal surfaces for which the resolution of 
the Calabi conjecture gives existence of an Einstein metric.
\end{rem}

In view of Theorem~\ref{t:main}, one can ask how many
smooth structures with Einstein metrics and how many
without, a given topological manifold has. On the one
hand, using for example the work of Fintushel--Stern, 
one can show that one has infinitely many choices for 
the smooth structures of the manifolds $Y_i$ in the proof 
of Theorem~\ref{t:main}, which remain distinct under blowing 
up points. Thus, one has infinitely many smooth manifolds one 
can use for each $X_i$, not admitting any Einstein metrics. 
On the other hand, it is known that there are homeomorphic
non-diffeomorphic minimal surfaces of general type,
cf~\cite{FM}, page 410, and the references cited there.
In fact, the number of distinct smooth structures among
sets of homeomorphic minimal surfaces of general type
can be arbitrarily large~\cite{Salve}. It is not hard to
check that all the examples in~\cite{Salve} and~\cite{FM}
have ample canonical bundle, and therefore have K\"ahler--Einstein
metrics of negative scalar curvature. However, all
those examples have $c_1^2 > 3\chi$, and can therefore
not be used as the $Z_i$ in the proof of Theorem~\ref{t:main}.
Still, those examples show that a given simply connected
topological manifold can have an arbitrarily large number
of smooth structures admitting Einstein metrics. Compare 
Theorem~\ref{t:top} below.

\section{Uniqueness for a given smooth structure}

We have seen that existence of Einstein metrics on closed
$4$--manifolds depends in an essential way on the smooth
structure. I believe that the issue of uniqueness, up
to the sign of the scalar curvature, is also
tied to the smooth structure. More specifically:

\begin{conj}
A closed smooth $4$--manifold admits Einstein metrics for 
at most one sign of the  scalar curvature.
\end{conj}

Such questions were raised in~\cite{einstein}, pages 18--19, 
and are also addressed in~\cite{catleb}. What is new here,
and in~\cite{catleb}, 
is that the answer depends on the smooth structure, and also 
seems to depend on the dimension. The conjecture is interesting 
because it is sharp --- it would be false if one did not fix the 
smooth structure, but only the underlying topological manifold:

\begin{thm}\label{t:top}
There are simply connected homeomorphic but non-diffeomor\-phic
smooth $4$--manifolds $X$ and $Y$, such that $X$
admits an Einstein metric of positive scalar curvature, and
$Y$ admits an Einstein metric of negative scalar curvature.
\end{thm}
\begin{proof}
We can take for $X$ the $8$--fold blowup of $\C P^2$. By the 
work of Tian--Yau~\cite{TY} this admits a K\"ahler--Einstein
metric of positive scalar curvature. 

For $Y$ we take the simply connected numerical Godeaux
surface constructed by Craighero--Gattazzo~\cite{CG} and
studied recently by Dolgachev--Werner~\cite{dolgwer}\footnote{I 
am grateful to Igor Dolgachev for telling me about this, and 
for providing an advance copy of~\cite{dolgwer}.}.
This has ample canonical bundle, so that by the work of Aubin
and Yau it admits a K\"ahler--Einstein metric of negative
scalar curvature. 

By Freedman's classification~\cite{freed}, $X$ and $Y$ are
homeomorphic. That they are not diffeomorphic is clear
from~\cite{invent}. The argument carried out there for the 
Barlow surface, cf~\cite{Barlow}, works even more easily for 
the Craighero--Gattazzo surface as there is no complication 
arising from $(-2)$--curves. Alternatively, the fact
that $X$ and $Y$ have K\"ahler--Einstein metrics of 
opposite signs implies via Seiberg--Witten theory
that they are non-diffeomorphic.
\end{proof}

In higher dimensions, homeomorphic non-diffeomorphic
manifolds with Einstein metrics are known~\cite{KS,WZ},
though in those examples all the metrics have positive
scalar curvature. Although the examples of~\cite{KS,WZ}
are consistent with a higher-dimensional analogue 
of the above conjecture, such a generalisation is false:

\begin{cor}
For every $i\geq 2$ there is a simply connected closed
$4i$--manifold which admits Einstein metrics of both
positive and negative scalar curvature.
\end{cor}
\begin{proof}
Let $X$ and $Y$ be as in Theorem~\ref{t:top}. Then the 
$i$--fold products $X^i=X\times\ldots\times X$ and 
$Y^i=Y\times\ldots\times Y$ have K\"ahler--Einstein metrics 
of positive, respectively negative, scalar curvature.
However, as $X$ and $Y$ are simply connected and 
homeomorphic, they are $h$--cobordant. Therefore, $X^i$
and $Y^i$ are also $h$--cobordant for all $i$, and for
$i\geq 2$ are diffeomorphic by the $h$--cobordism
theorem.
\end{proof}

Theorem~\ref{t:top} and Corollary 2.3 have also been 
proved independently by Catanese and LeBrun~\cite{catleb}.
Instead of the Craighero--Gattazzo surface they use the 
Barlow surface, showing that it has deformations with
ample canonical bundle. Conjecturally, the Barlow and 
Craighero--Gattazzo surfaces are deformation equivalent,
and therefore diffeomorphic.

\medskip

{\sl Acknowledgement:}\stdspace I am grateful to the Department of
Mathematics at Brown University and to the 
Max--Planck--Institut f\"ur Mathematik in Bonn
 for hospitality and support.

%
%

%
%

\begin{thebibliography}


\bibitem{Barlow}
{\bf R\ Barlow}, {\em A simply connected surface of general type with
$p_g = 0$}, Invent.~Math.  79 (1985) 293--301

\bibitem{BPV}
{\bf W\ Barth}, {\bf C\ Peters}, {\bf A\ Van~de~Ven}, {\sl Compact Complex
Surfaces}, Springer--Verlag (1984)

\bibitem{einstein}
{\bf A\,L\ Besse}, {\sl Einstein Manifolds}, Springer--Verlag (1987)

\bibitem{catleb}
{\bf F\ Catanese}, {\bf C\ LeBrun}, {\em On the scalar curvature of Einstein
manifolds}, preprint (May 1997)

\bibitem{CG} 
{\bf P\ Craighero}, {\bf R\ Gattazzo}, {\em Quintic surfaces of $\C
P^3$ having a nonsingular model with $q=p_g=0$, $P_2 \neq 0$,}
Rend.~Sem.~Mat.~Univ.~Padova {91} (1994) 187--198

\bibitem{dolgwer}
{\bf I\ Dolgachev}, {\bf C\ Werner}, {\em A simply connected numerical
Godeaux surface with ample canonical class}, preprint (April 1997)

\bibitem{bams}
{\bf S\,K\ Donaldson}, {\em The Seiberg--Witten equations and $4$--manifold
topology}, Bull.~Amer.~Math.~Soc. {33} (1996) 45--70

\bibitem{FS}
{\bf R\ Fintushel}, {\bf R\,J\ Stern}, {\em Immersed spheres in $4$--manifolds
and the immersed Thom conjecture}, Proc. of G\"okova Geometry--Topology
Conference 1994, Turkish J. of Math. {19} (2) (1995) 27--39

\bibitem{freed}
{\bf M\,H\ Freedman}, {\em The topology of four--manifolds}, J.~Differential
Geometry {17} (1982) 357--454

\bibitem{FM}
{\bf R\ Friedman}, {\bf J\,W\ Morgan}, {\sl Smooth Four--Manifolds and 
Complex Surfaces}, Springer--Verlag (1994)

\bibitem{gromov}
{\bf M~Gromov}, {\em Volume and bounded cohomology}, 
Publ.~Math.~I.H.E.S. {56} (1982) 5--99

\bibitem{HT}
{\bf N\,J~Hitchin}, {\em Compact four--dimensional Einstein manifolds}, 
J.~Differential Geometry {9} (1974) 435--441

\bibitem{invent}
{\bf D~Kotschick}, {\em On manifolds homeomorphic to $\C P^2\# 8\blowup$}, 
Invent.~math. {95} (1989) 591--600

\bibitem{orient}
{\bf D~Kotschick}, {\em Orientation-reversing homeomorphisms in surface
geography}, Math. Annalen {292} (1992) 375--381

\bibitem{KKM}
{\bf D~Kotschick}, {\bf P\,B Kronheimer}, {\bf T\,S~Mrowka}, monograph in 
preparation

\bibitem{KMT}
{\bf D~Kotschick}, {\bf J\,W~Morgan}, {\bf C\,H~Taubes}, {\em Four--manifolds 
without symplectic structures but with non-trivial Seiberg--Witten 
invariants}, Math. Research Letters {2} (1995) 119--124

\bibitem{KS}
{\bf M~Kreck}, {\bf S~Stolz}, {\em A diffeomorphism classification of
$7$--dimensional homogeneous Einstein manifolds with
$SU(3)\times SU(2)\times U(1)$ symmetry}, Ann.~of Math. {127}
(1988) 373--388

\bibitem{lebrun}
{\bf C~LeBrun}, {\em Four--manifolds without Einstein metrics},
Math. Research Letters {3} (1996) 133--147

\bibitem{per}
{\bf U~Persson}, {\em Chern invariants of surfaces of general type},
Comp.~Math. {43} (1981) 3--58

\bibitem{Salve}
{\bf M~Salvetti}, {\em On the number of non-equivalent differentiable
structures on $4$--manifolds}, Manuscr.~Math. {63} (1989)
157--171

\bibitem{BCG}
{\bf A~Sambusetti}, {\em An obstruction to the existence of Einstein
metrics on $4$--manifolds}, C.~R.~Acad.~Sci.~Paris {322}
(1996) 1213--1218

\bibitem{TY}
{\bf G~Tian}, {\bf S-T~Yau}, {\em K\"ahler--Einstein metrics on complex
surfaces  with $c_1 >0$}, Comm.~Math.~Phys. {112} (1987)
175--203

\bibitem{WZ}
{\bf McK\,Y~Wang}, {\bf W~Ziller}, {\em Einstein metrics on principal
torus bundles}, J.~Differential Geometry {31} (1990) 215--248

\bibitem{witten}
{\bf E~Witten}, {\em Monopoles and four--manifolds}, Math. Research Letters
{1} (1994) 769--796

\end{thebibliography}
\end{document}